\documentclass[reqno]{amsart}
\usepackage{amssymb,url}
\usepackage{hyperref}
\date{\bf April 18, 2008}
\theoremstyle{plain}
\newtheorem{theorem}{Theorem}
\newtheorem{lemma}{Lemma}
\newtheorem{definition}{Definition}
\newtheorem{corollary}{Corollary}
\newtheorem{proposition}{Proposition}
\theoremstyle{remark}
\newtheorem{remark}{Remark}

\allowdisplaybreaks
\DeclareMathOperator{\T}{\mathbb T}
\DeclareMathOperator{\sigt}{\sigma (t)}
\DeclareMathOperator{\rt}{\rho (t)}

\DeclareMathOperator{\fd}{{\it f}^\Delta (t)}

\DeclareMathOperator{\fs}{{\it f}^\sigma (t)}

\DeclareMathOperator{\Fd}{{\it F}^\Delta (t)}

\begin{document}

\title{A generalization of Ostrowski inequality  on time scales for $k$ points}

\author[W. J. Liu]{Wenjun Liu}
\address[W. J. Liu]{College of Mathematics and Physics\\
Nanjing University of Information Science and Technology \\
Nanjing 210044, China} \email{\href{mailto: W. J. Liu
<wjliu@nuist.edu.cn>}{wjliu@nuist.edu.cn}}

\author[Q. A. Ng\^{o}]{Qu\^{o}\hspace{-0.5ex}\llap{\raise 1ex\hbox{\'{}}}\hspace{0.5ex}c Anh Ng\^{o}}
\address[Q. A. Ng\^{o}]{Department of Mathematics, Mechanics and Informatics\\
College of Science\\ Vi\d{\^{e}}t Nam National University\\ H\`{a}
N\d{\^{o}}i, Vi\d{\^{e}}t Nam} \email{\href{mailto: Q. A. Ng\^{o}
<bookworm\_vn@yahoo.com>}{bookworm\_vn@yahoo.com}}



\begin{abstract}
In this paper we first generalize the Ostrowski inequality on time
scales for $k$ points and then unify corresponding continuous and
discrete versions. We also point out some particular Ostrowski type
inequalities on time scales as special cases.
\end{abstract}

\thanks{This paper was typeset using \AmS-\LaTeX}

\maketitle

\section{introduction}

In 1938, A. Ostrowski proved the following interesting integral
inequality  which has received considerable attention from many
researchers \cite{l1, l3, l4, mpf, mpf2}.
 \begin{theorem}\label{t1}
 Let $f: [a,b]\rightarrow \mathbb{R}$ be continuous on $[a,b]$ and differentiable in $(a,b)$ and its derivative $f': (a,b)\rightarrow \mathbb{R}$ is bounded in $(a,b)$, that is, $\|f'\|_\infty:=\sup\limits_{t\in(a,b)}|f'(x)|< \infty$. Then for any $x\in [a,b]$, we have the inequality:
\[
\left|\int\limits_a^bf(t)dt-f(x)(b-a)\right|\leq
\left(\frac{(b-a)^2}{4}+ \left(x-\frac{a+b}{2}\right)^2
\right)\|f'\|_\infty.\label{1.1}
 \]
The inequality is sharp in the sense that the constant $\frac{1}{4}$ cannot be replaced by a smaller one.
\end{theorem}

The development of the theory of time scales was initiated by Hilger
\cite{h1988} in 1988 as a theory capable to contain both difference
and differential calculus in a consistent way. Since then, many
authors have studied the theory of certain integral inequalities or
dynamic equations on time scales. For example, we refer the reader
to \cite{abp2001, bm2007, bm2008, gz, ln, OSY, sz, wyy}. In
\cite{bm2008}, Bohner and Matthews established  the following
so-called Ostrowski inequality on time scales.

\begin{theorem}[See \cite{bm2008}, Theorem 3.5]\label{th2}
Let $a, b, x, t\in \T$, $a<b$ and $f: [a, b]\rightarrow \mathbb{R}$ be differentiable. Then
\begin{equation}\label{eq1}
 \left|\int\limits_a^b f^\sigma(t)\Delta t-f(x)(b-a)\right|  \leq M\Big(h_2(x,a)+h_2(x,b)\Big),
\end{equation}
where $h_2(\cdot, \cdot)$ is defined by Definition \ref{de7} below and $M=\sup_{a<x<b}|f^\Delta(x)|.$ This inequality is sharp in the
sense that the right-hand side of (\ref{eq1}) cannot be replaced by a smaller one.
\end{theorem}

In the present paper we shall first  generalize the above Ostrowski
inequality on time scales for $k$ points $x_1, x_2, \cdots, x_k$ and
then unify corresponding continuous and discrete versions. We also
point out some particular Ostrowski type inequalities on time scales
as special cases.

\section{Time scales essentials}

Now we briefly introduce the time scales theory and refer the reader to Hilger \cite{h1988} and the books \cite{bp2001, bp2003, LSK} for
further details.

\begin{definition}
{\it A time scale} $\T$ is an arbitrary nonempty closed subset of real numbers.
\end{definition}


\begin{definition}
 For $t \in \T$, we define the {\it forward jump operator} $\sigma : \T \to \T$ by
$
\sigt = \inf \left\{ {s \in \T:s > t} \right\},
$
while the {\it backward jump operator} $\rho : \T \to \T$ is defined by
$
\rt = \sup \left\{ {s \in \T:s < t} \right\}.
$
If $\sigt >t $, then we say that $t$ is {\it right-scattered}, while if $\rt < t$ then we say that $t$ is {\it left-scattered}.
\end{definition}

Points that are right-scattered and left-scattered at the same time are called isolated. If $\sigt = t$, the $t$ is called {\it
right-dense}, and if $\rt = t$ then $t$ is called {\it left-dense}. Points that are both right-dense and left-dense are called dense.

\begin{definition}
Let $t \in \T$, then two mappings $\mu ,\nu :\T \to \left( {0, +
\infty } \right)$ satisfying $$ \mu \left( t \right): = \sigt - t,
\, \nu \left( t \right): = t - \rt$$ are called the {\it graininess
functions}.
\end{definition}

We now introduce the set $\T^\kappa$ which is derived from the time scales $\T$ as follows. If $\T$ has a left-scattered maximum $t$,
then $\T^\kappa := \T -\{t\}$, otherwise $\T^\kappa := \T$. Furthermore for a function $f : \T \to \mathbb R$, we define the function $f^\sigma : \T \to \mathbb R$ by $\fs = f(\sigma(t))$ for all $t \in \T$.

\begin{definition}
Let $f : \T \to \mathbb R$ be a function on time scales. Then for $t \in \T^\kappa$, we define $\fd$ to be the number, if one exists, such that for all $\varepsilon >0$ there is a neighborhood $U$ of $t$ such that for all $s \in U$
\[
\left| {\fs - f\left( s \right) - \fd \left( {\sigt - s} \right)}
\right| \leq \varepsilon \left| {\sigt - s} \right|.
\]
We say that $f$ is $\Delta$-differentiable on $\T^\kappa$ provided $\fd$ exists for all $t \in \T^\kappa$.
\end{definition}


\begin{definition} A mapping $f : \T \to \mathbb R$ is called {\it rd-continuous} (denoted by $C_{rd}$) provided if it satisfies
  \begin{enumerate}
   \item $f$ is continuous at each right-dense point or maximal element of $\T$.
   \item The left-sided limit $\mathop {\lim }\limits_{s \to t - } f\left( s \right) = f\left( {t - } \right)$ exists at each left-dense point $t$ of $\T$.
  \end{enumerate}
\end{definition}

\begin{remark}
It follows from Theorem 1.74 of Bohner and Peterson \cite{bp2001} that every rd-continuous function has an anti-derivative.
\end{remark}

\begin{definition} A function $F : \T \to \mathbb R$ is called a $\Delta$-antiderivative of $f : \T \to \mathbb R$ provided $\Fd =f(t)$ holds for all $t \in \T^\kappa$. Then the $\Delta$-integral of $f$ is defined by
\[
\int\limits_a^b {f\left( t \right)\Delta t} = F\left( b \right) - F\left( a \right).
\]
\end{definition}

\begin{proposition} \label{pro1} Let $f, g$ be rd-continuous, $a, b, c\in \T$ and $\alpha, \beta\in \mathbb{R}$. Then
\begin{enumerate}
  \item $\int\limits_a^b {\Big(\alpha f(t)+\beta g(t)\Big)\Delta t} = \alpha\int\limits_a^b {f(t)\Delta t}+\beta\int\limits_a^b{g(t)\Delta t}$,
  \item $\int\limits_a^b{f(t)\Delta t}=-\int\limits_b^a {f(t)\Delta t}$,
  \item $\int\limits_a^b{f(t)\Delta t}=\int\limits_a^c {f(t)\Delta t}+\int\limits_c^b{f(t)\Delta t}$,
  \item $\int\limits_a^b{f(t)g^\Delta(t)\Delta t}=(fg)(b)-(fg)(a)-\int\limits_a^b{f^\Delta(t)g(\sigma(t))\Delta t},$
  \item $\int\limits_a^a{f(t)\Delta t}=0$.
 \end{enumerate}
\end{proposition}

\begin{definition}\label{de7}
Let $h_k : \T^2 \to \mathbb R$, $k \in \mathbb N_0$ be defined by
\[
h_0 \left( {t,s} \right) = 1 \quad {\text{ for all }} \quad s,t \in
\T
\]
and then recursively by
\[
h_{k + 1} \left( {t,s} \right) = \int\limits_s^t {h_k \left( {\tau
,s} \right)\Delta \tau } \quad {\text{ for all }} \quad s,t \in \T.
\]
\end{definition}

\section{The generalized Ostrowski inequality on time scales}

Throughout this section, we suppose that $\T$ is a time scale and an interval means the intersection of real interval with the given time scale. We are in a position to state our main result.

\begin{theorem}\label{th3}
Suppose that
\begin{enumerate}
 \item $a, b\in \T$, $I_k: a=x_0<x_1<\cdots<x_{k-1}<x_k=b$ is a division of the interval $[a, b]$ for $x_0, x_1, \dots, x_k\in \T$;
 \item $\alpha_i \in \T$ $(i=0, \dots, k+1)$ is "$k+2$" points so that $\alpha_0=a$, $\alpha_i\in [x_{i-1}, x_i]$ $(i=1, \dots, k)$ and $\alpha_{k+1}=b$;
 \item $f: [a, b]\rightarrow \mathbb{R}$ is differentiable.
 \end{enumerate}
Then we have
\begin{equation}\label{eq2}
 \left|\int\limits_a^b f^\sigma(t)\Delta
 t-\sum_{i=0}^k(\alpha_{i+1}-\alpha_i)f(x_i)\right|
  \leq M\sum_{i=0}^{k-1}\Big(h_2(x_{i}, \alpha_{i+1})+h_2(x_{i+1}, \alpha_{i+1})\Big),
\end{equation}
where $$M=\sup_{a<x<b}|f^\Delta(x)|.$$ This inequality is sharp in the sense that the right-hand side of (\ref{eq2}) cannot be replaced by a smaller one.
\end{theorem}

To prove Theorem \ref{th3}, we need the following Generalized Montgomery Identity.

\begin{lemma}[Generalized Montgomery Identity]\label{le1}
Under the assumptions of Theorem \ref{th3}, we have
\begin{equation}
\sum_{i=0}^k(\alpha_{i+1}-\alpha_i)f(x_i)=\int\limits_a^b
f^\sigma(t)\Delta
 t+ \int\limits_a^b K(t,I_k)f^\Delta(t)\Delta
 t,\label{eq3}
\end{equation}
where
\begin{equation}
 K(t,I_k)=\left\{ \begin{array}{cl}
 t-\alpha_1, &t\in [a, x_1), \hfill \\
  t-\alpha_2, &t\in [x_1, x_2), \hfill \\
  \cdots & \cdots\hfill \\
  t-\alpha_{k-1}, &t\in [x_{k-2}, x_{k-1}), \hfill \\
  t-\alpha_k, &t\in [x_{k-1}, b). \hfill
\end{array} \right.\label{eq4}
\end{equation}
\end{lemma}

\begin{proof} Integrating by parts and applying Proposition \ref{pro1}, we have
\begin{align*}
 \int\limits_a^b K(t,I_k)f^\Delta(t)&\Delta t=\sum_{i=0}^{k-1}\int\limits_{x_i}^{x_{i+1}} K(t,I_k)f^\Delta(t)\Delta  t\\
 = &\sum_{i=0}^{k-1}\int\limits_{x_i}^{x_{i+1}} (t-\alpha_{i+1})f^\Delta(t)\Delta
 t\\
 = &\sum_{i=0}^{k-1}\left( (x_{i+1}-\alpha_{i+1})f(x_{i+1}) - (x_i - \alpha_{i+1})f(x_i)
 -\int\limits_{x_i}^{x_{i+1}} f^\sigma(t)\Delta t\right) \\
 = &\sum_{i=0}^{k-1}\left((\alpha_{i+1}-x_i)f(x_i)+(x_{i+1}-\alpha_{i+1})f(x_{i+1})
 -\int\limits_{x_i}^{x_{i+1}} f^\sigma(t)\Delta t\right) \\
 = &(\alpha_1-a)f(a)+\sum_{i=1}^{k-1}(\alpha_{i+1}-x_i)f(x_i)+ \sum_{i=0}^{k-2}(x_{i+1}-\alpha_{i+1})f(x_{i+1})\\
 & \quad\quad +(b-\alpha_k)f(b)-\int\limits_{a}^{b} f^\sigma(t)\Delta
 t\\
 = &(\alpha_1-a)f(a)+\sum_{i=1}^{k-1}(\alpha_{i+1}-\alpha_i)f(x_i)+(b-\alpha_k)f(b)-\int\limits_{a}^{b} f^\sigma(t)\Delta
 t\\
 = & \sum_{i=0}^{k}(\alpha_{i+1}-\alpha_i)f(x_i) -\int\limits_{a}^{b} f^\sigma(t)\Delta
 t,
\end{align*}
i.e., (\ref{eq3}) holds.
\end{proof}

\begin{proof}[Proof of Theorem~\ref{th3}]
By applying Lemma \ref{le1}, we get
\begin{align*}
\Bigg|\int\limits_a^b f^\sigma(t)\Delta t-\sum_{i=0}^k&(\alpha_{i+1}-\alpha_i)f(x_i)\Bigg|\\
= &\left|\int\limits_a^b K(t,I_k)f^\Delta(t)\Delta
 t\right|
 = \left|\sum_{i=0}^{k-1}\int\limits_{x_i}^{x_{i+1}} K(t,I_k)f^\Delta(t)\Delta
 t\right|\\
 \leq &\sum_{i=0}^{k-1}\int\limits_{x_i}^{x_{i+1}} |K(t,I_k)|\left|f^\Delta(t)\right|\Delta
 t\leq M\sum_{i=0}^{k-1}\int\limits_{x_i}^{x_{i+1}} |t-\alpha_{i+1}| \Delta
 t\\
 = &M\sum_{i=0}^{k-1}\left(\int\limits_{x_i}^{\alpha_{i+1}} (\alpha_{i+1}-t) \Delta
 t+\int\limits_{\alpha_{i+1}}^{x_{i+1}} (t-\alpha_{i+1}) \Delta
 t\right)\\
 = &M\sum_{i=0}^{k-1} \Big(h_2(x_{i}, \alpha_{i+1})+h_2(x_{i+1}, \alpha_{i+1}) \Big).
\end{align*}

To prove the sharpness of this inequality, let $f(t)=t,$ $x_0=a,$ $x_1=b,$ $\alpha_0=a,$ $\alpha_1=b$,  $\alpha_2=b$. It follows that $M=1$. Starting with the left-hand  side of (\ref{eq2}), we have
 \begin{align*}
 \left|\int\limits_a^b f^\sigma(t)\Delta
 t-\sum_{i=0}^k(\alpha_{i+1}-\alpha_i)f(x_i)\right| = &\left|\int\limits_a^b \sigma(t)\Delta
 t- \Big((b-a)a+(b-b)b \Big)\right|
 \\
 = &\left|\int\limits_a^b (\sigma(t)+t)\Delta
 t-\int\limits_a^b t\Delta
 t-(b-a)a \right| \\
 = &\left|\int\limits_a^b (t^2)^\Delta\Delta
 t-\int\limits_a^b t\Delta
 t-(b-a)a \right|\\
 = &\left| (b-a)a -\int\limits_a^b t\Delta
 t\right|.
\end{align*}
Starting with the right-hand side of (\ref{eq2}), we have
 \begin{align*}
 M\sum_{i=0}^{k-1}\left(h_2( \alpha_{i+1}, x_i)+h_2( \alpha_{i+1},
 x_{i+1})\right)
 = &h_2(x_0, \alpha_1)+h_2(x_1, \alpha_1) \\
 = &h_2(a, b)+h_2(b, b)\\
 = &\int_b^a (t-b)\Delta t+\int_b^b (t-b)\Delta t\\
 = &\int_b^a  t \Delta t - \int_b^a b\Delta t\\
 = &b(b-a)-\int_a^b t \Delta t.
\end{align*}
Therefore in this particular case
$$\left|\int\limits_a^b f^\sigma(t)\Delta
 t-\sum_{i=0}^k(\alpha_{i+1}-\alpha_i)f(x_i)\right|
  \geq M\sum_{i=0}^{k-1}\Big(h_2(\alpha_{i+1},x_{i})+h_2(\alpha_{i+1},x_{i+1})\Big)$$
and by (\ref{eq2}) also $$\left|\int\limits_a^b f^\sigma(t)\Delta
 t-\sum_{i=0}^k(\alpha_{i+1}-\alpha_i)f(x_i)\right|
  \leq
  M\sum_{i=0}^{k-1}\Big(h_2(\alpha_{i+1},x_{i})+h_2(\alpha_{i+1},x_{i+1})\Big).$$
So the sharpness of the inequality (\ref{eq2}) is shown.
\end{proof}

If we apply the inequality (\ref{eq2}) to different time scales, we will get some well-known and some new results.

\begin{corollary}[Continuous case]\label{co1}
Let $\mathbb{T }= \mathbb{R}$. Then our delta integral is the usual Riemann integral from calculus. Hence,
\[
h_2 \left( {t,s} \right) = \frac{{\left( {t - s} \right)^2 }}{2}, \quad {\text{ for all }} \quad t, s \in \mathbb R.
\]
This leads us to state the following inequality
\begin{equation}\label{eq5}
\begin{split}
 \Bigg|\int\limits_a^b f(t)\Delta t-\sum_{i=0}^k(\alpha_{i+1} & -\alpha_i)f(x_i)\Bigg| \\
  \leq &M\left(\frac{1}{4}\sum_{i=0}^{k-1}(x_{i+1}-x_i)^2+\sum_{i=0}^{k-1}\left(\alpha_{i+1}-\frac{x_i+x_{i+1}}{2}\right)^2\right),
\end{split}
\end{equation}
where $M=\sup\limits_{a<x<b}|f'(x)| $ and the constant $\frac{1}{4}$ in the right-hand side is the best possible.
\end{corollary}

\begin{remark}
The inequality (\ref{eq5}) is exactly the generalized Ostrowski
inequality shown in \cite{d}.
\end{remark}

\begin{corollary}[Discrete case]\label{co2}
Let $\mathbb{T }= \mathbb{Z}$, $a=0$, $b=n$. Suppose that
\begin{enumerate}
   \item $I_k:
   0=j_0<j_1<\cdots<j_{k-1}<j_k=n$ is a division of
   $[0, n]\cap\mathbb{Z}$ for $j_0, k_1, \dots, j_k\in \mathbb{Z};$
   \item $p_i \in \mathbb{Z}$ $(i=0, \dots, k+1)$ is "$k+2$"
   points so that $p_0=0$, $p_i\in [j_{i-1}, j_i]\cap\mathbb{Z}$ $(i=1, \dots,
   k)$ and $p_{k+1}=n$;
\item $f(k)=x_k$.
\end{enumerate}
Then, we have
\begin{align*} &\left|\sum_{j=1}^nx_j-\sum_{i=0}^k(p_{i+1}-p_i)x_{j_i}\right|\\ \leq
&M\left(\frac{1}{4}\sum_{i=0}^{k-1}(j_{i+1}-j_i)^2+\sum_{i=0}^{k-1}\left(p_{i+1}-\frac{j_i+j_{i+1}}{2}\right)^2
+\sum_{i=0}^{k-1}\left(p_{i+1}-\frac{j_i+j_{i+1}}{2}\right)\right)
\end{align*}
for all $i=\overline{1,n}$, where $M=\sup\limits_{i=1,\cdots, n-1} \left|\Delta x_i\right|$ and the constant $\frac{1}{4}$ in the right-hand side is the best possible.
\end{corollary}

\begin{proof}
 It is known that
\[
h_k \left( {t,s} \right) = \left( {\begin{array}{*{20}c}
  {t - s} \\
  k \\
 \end{array} } \right) , \quad {\text{ for all }} \quad t,s \in \mathbb Z.
\]
Therefore,
\[
h_2 \left( {j_i, p_{i+1}} \right) = \left( {\begin{array}{*{20}c}
  j_i- p_{i+1} \\
  2 \\
 \end{array} } \right) = \frac{{(j_i- p_{i+1})(j_i- p_{i+1}-1)}}
{2}
\]
and
\[
h_2 \left( {j_{i+1}, p_{i+1}} \right) = \left(
{\begin{array}{*{20}c}
  j_{i+1}- p_{i+1} \\
  2 \\
 \end{array} } \right) = \frac{{(j_{i+1}- p_{i+1})(j_{i+1}- p_{i+1}-1)}}
{2}.
\]
The conclusion is obtained by some easy calculation.
\end{proof}

\begin{corollary}[Quantum calculus case]\label{co3}
Let $\mathbb{T }= q^{\mathbb{N}_0}$, $q>1$, $a=q^m, b=q^n$ with $m<n$. Suppose that
\begin{enumerate}
   \item $I_k:
   q^m=q^{j_0}<q^{j_1}<\cdots<q^{j_{k-1}}<q^{j_k}=q^n$ is a division of
   $[q^m, q^n]\cap q^{\mathbb{N}_0}$ for $j_0, k_1, \dots, j_k\in {\mathbb{N}_0};$
   \item $q^{p_i} \in q^{\mathbb{N}_0}$ $(i=0, \dots, k+1)$ is "$k+2$"
   points so that $q^{p_0}=q^m$, $q^{p_i}\in [q^{j_{i-1}}, q^{j_i}]\cap q^{\mathbb{N}_0}$ $(i=1, \dots,
   k)$ and $q^{p_{k+1}}=q^m$;
\item $f: [q^m, q^n]\rightarrow \mathbb{R}$ is differentiable.
\end{enumerate}
Then, we have
\begin{multline*}
 \Bigg|\int\limits_{q^m}^{q^n} f^\sigma(t)\Delta t -\sum_{i=0}^k(q^{p_{i+1}}-q^{p_i})f\left(q^{j_i}\right)
 \Bigg| \\
\leq \frac{2M}{1+q}\sum_{i=0}^{k-1}\Bigg(\left(q^{j_i}-\frac{\frac{1+q}{2}\left(q^{p_i}+q^{p_{i+1}}\right)}{2}\right)^2 + \\
\left. \frac{2\left(q^{2p_i}+q^{2p_{i+1}}\right)-(\frac{1+q}{2})^2\left(q^{p_i}+q^{p_{i+1}}\right)^2}{4}+q^{2j_i}(q-1)\right),
\end{multline*}
where
$$M=\sup\limits_{q^m<t<q^n}\left|\frac{f(qt)-f(t)}{(q-1)(t)}\right|$$
and the constant $\frac{1}{4}$ in the right-hand side is the best possible.
\end{corollary}
\begin{proof}
 In this situation, one has
\[
h_k \left( {t,s} \right) = \prod\limits_{\nu = 0}^{k - 1} {\frac{{t
- q^\nu s}} {{\sum\limits_{\mu = 0}^\nu {q^\mu } }}}, \quad
{\text{ for all }} \quad t,s \in q^{\mathbb{N}_0}.
\]
Therefore,
\[
h_2 \left( {q^{j_i}, q^{p_{i+1}}} \right) = \frac{{\left( {q^{j_i}-
q^{p_{i+1}}} \right)\left( {q^{j_i}- q^{p_{i+1}+1}} \right)}} {{1 +
q}}
\]
and
\[
h_2 \left( {q^{j_{i+1}}, q^{p_{i+1}}} \right) = \frac{{\left(
{q^{j_{i+1}}- q^{p_{i+1}}} \right)\left( {q^{j_{i+1}}-
q^{p_{i+1}+1}} \right)}} {{1 + q}}.
\]
The conclusion is easy obtained by some simple calculation.
\end{proof}

\section{Some particular Ostrowski type inequalities on time scales}

In this section we point out some particular Ostrowski type inequalities on time scales as special cases, such as: {\it rectangle inequality}
on time scales, {\it trapezoid inequality} on time scales, {\it mid-point inequality} on time scales, {\it Simpson inequality} on time scales,
 {\it averaged mid-point-trapezoid inequality} on time scales and
others.

Throughout this section, we always assume $\T$ is a time scale; $a, b \in \T$ with $a < b$; $f: [a, b]\rightarrow \mathbb{R}$ is differentiable.
 We denote $$M=\sup_{a<x<b}|f^\Delta(x)|.$$

\begin{proposition}\label{pro2}
Suppose that $\alpha \in [a,b] \cap \T$. Then we have the sharp rectangle inequality on time scales
\begin{equation}\label{eq6}
 \left|\int\limits_a^b f^\sigma(t)\Delta
 t- \Big((\alpha-a)f(a)+(b-\alpha)f(b) \Big) \right| \leq M \Big(h_2(a, \alpha)+h_2(b, \alpha) \Big).
\end{equation}
\end{proposition}

\begin{proof}
We choose $k=1$, $x_0=a$, $x_1=b$, $\alpha_0=a$, $\alpha_1=\alpha$ and $\alpha_2=b$ in Theorem \ref{th3} to get the result.
\end{proof}

\begin{remark}

\begin{enumerate}
  \item [(a)] If we choose $\alpha=b$ in (\ref{eq6}), we get the sharp left rectangle inequality on time scales
\begin{equation}\label{eq7}
 \left|\int\limits_a^b f^\sigma(t)\Delta
 t-(b-a)f(a)\right|
  \leq M h_2(a, b).
\end{equation}

  \item [(b)] If we choose $\alpha=a$ in (\ref{eq6}), we get the sharp right rectangle inequality on time scales
\begin{equation}\label{eq8}
 \left|\int\limits_a^b f^\sigma(t)\Delta
 t-(b-a)f(b)\right|
  \leq M h_2(a, b).
\end{equation}

  \item [(c)] If we choose $\alpha=\frac{a+b}{2}$ in (\ref{eq6}), we get the sharp trapezoid inequality on time scales
\begin{equation}\label{eq9}
 \left|\int\limits_a^b f^\sigma(t)\Delta
 t-\frac{f(a)+f(b)}{2}(b-a)\right|
  \leq M\left(h_2\left(a, \frac{a+b}{2}\right)+h_2\left(b, \frac{a+b}{2}\right)\right).
\end{equation}
\end{enumerate}
\end{remark}

\begin{proposition}\label{pro3}
Suppose that $x\in [a, b] \cap \T$, $\alpha_1\in [a, x] \cap \T$, $\alpha_2\in [x, b] \cap \T$. Then we have the sharp inequality on time scales
\begin{equation}\label{eq10}
\begin{split}
 \Bigg|\int\limits_a^b f^\sigma(t)\Delta t -\Big((\alpha_1-&a)f(a)+(\alpha_2-\alpha_1)f(x)+(b-\alpha_2)f(b)\Big)\Bigg| \\
  \leq &M\Big(h_2(a, \alpha_1)+h_2(x, \alpha_1)+h_2(x, \alpha_2)+h_2(b, \alpha_2)\Big).
\end{split}
\end{equation}
\end{proposition}

\begin{proof}
We choose $k=2$, $x_0=a$, $x_1=x$, $x_2=b$ and $\alpha_i$ ($i = \overline{0, 3}$) is as in Theorem \ref{th3} to get the result.
\end{proof}

\begin{remark}
\begin{enumerate}
  \item[(a)] If we choose $\alpha_1=a$ and $\alpha_2=b$ in Proposition \ref{pro3}, we get exactly Theorem \ref{th2}. Therefore, Theorem \ref{th3} is a generalization of Theorem 3.5 in \cite{bm2008}.
  \item[(b)] If we choose $x=\frac{a+b}{2}$ in (\ref{eq1}), we get the sharp mid-point inequality on time scales
\begin{equation}\label{eq11}
 \left|\int\limits_a^b f^\sigma(t)\Delta  t-f\left(\frac{a+b}{2}\right)(b-a)\right|   \leq M\left(h_2\left(\frac{a+b}{2}, a\right)+h_2\left(\frac{a+b}{2}, b\right)\right).
\end{equation}
\end{enumerate}
\end{remark}

\begin{corollary}\label{co4}
Suppose that $\alpha_1=\frac{5a+b}{6}\in \T$, $\alpha_2=\frac{a+5b}{6}\in \T$, and $x\in \left[\frac{5a+b}{6},\frac{a+5b}{6}\right] \cap \T$. Then we have the sharp inequality on time scales
\begin{equation}\label{eq12}
\begin{split}
 &\left|\int\limits_a^b f^\sigma(t)\Delta
 t-\frac{b-a}{3}\left(\frac{f(a)+f(b)}{2}+2f(x)\right)\right| \\
  \leq &M\left(h_2\left(a, \frac{5a+b}{6}\right)+h_2\left(x, \frac{5a+b}{6}\right)
  +h_2\left(x, \frac{a+5b}{6}\right)+h_2\left(b, \frac{a+5b}{6}\right)\right).
\end{split}
\end{equation}
\end{corollary}

\begin{remark}
If we choose $x=\frac{a+b}{2}$ in (\ref{eq12}), we get the sharp Simpson inequality on time scales
\begin{align*}
 &\left|\int\limits_a^b f^\sigma(t)\Delta
 t-\frac{b-a}{3}\left(\frac{f(a)+f(b)}{2}+2f\left(\frac{a+b}{2}\right)\right)\right| \\
  \leq &M\left(h_2\left(a, \frac{5a+b}{6}\right)+h_2\left(\frac{a+b}{2}, \frac{5a+b}{6}\right)
  +h_2\left(\frac{a+b}{2}, \frac{a+5b}{6}\right)+h_2\left(b,
  \frac{a+5b}{6}\right)\right).
\end{align*}
\end{remark}

\begin{corollary}\label{co5}
Suppose that $\alpha_1\in \left[a, \frac{a+b}{2}\right] \cap \T$ and $\alpha_2\in \left[\frac{a+b}{2}, b\right] \cap \T$. Then we have the sharp inequality on time scales
\begin{equation}\label{eq13}
\begin{split}
 &\left|\int\limits_a^b f^\sigma(t)\Delta
 t-\left((\alpha_1-a)f(a)+(\alpha_2-\alpha_1)f\left(\frac{a+b}{2}\right)+(b-\alpha_2)f(b)\right)\right| \\
  \leq &M\left(h_2(a, \alpha_1)+h_2\left(\frac{a+b}{2}, \alpha_1\right)+h_2\left(\frac{a+b}{2},
  \alpha_2\right)+h_2(b, \alpha_2)\right).
\end{split}
\end{equation}
\end{corollary}

\begin{remark}
If we choose $\alpha_1=\frac{3a+b}{4}$ and $\alpha_2=\frac{a+3b}{4}$ in (\ref{eq13}), we get the sharp averaged mid-point-trapezoid inequality on time scales
\[
\begin{split}
  &\left|\int\limits_a^b f^\sigma(t)\Delta  t-\frac{b-a}{2}\left(\frac{f(a)+f(b)}{2}+f\left(\frac{a+b}{2}\right)\right)\right| \\
 \leq &M\left(h_2\left(a, \frac{3a+b}{4}\right)+h_2\left(\frac{a+b}{2}, \frac{3a+b}{4}\right)+h_2\left(\frac{a+b}{2},   \frac{a+3b}{4}\right)+h_2\left(b, \frac{a+3b}{4}\right)\right).
\end{split}
\]
\end{remark}

\medskip
\section*{Acknowledgements}

This work was supported by the Science Research Foundation of Nanjing University of Information Science and Technology and the
Natural Science Foundation of Jiangsu Province Education Department under Grant No.07KJD510133.

\end{document}